\documentclass[11pt]{amsart}

\usepackage{amssymb,amsmath,graphicx,caption,wrapfig}

\usepackage{color}

\newcommand\bigfrown[2][\textstyle]{\ensuremath{%
 \array[b]{c}\text{\resizebox{3ex}{.7ex}{$#1\frown$}}\\[-1.2ex]#1#2\endarray}}

\def\eq{\mathrel{\Leftrightarrow}}

\title[Dividing the Circle]{Dividing the Circle}

\begin{document}

\author[Hugo Tavares]{Hugo Tavares}\thanks{}
\address{\footnotesize{Hugo Tavares \newline \indent  CAMGSD, Mathematics Department \newline \indent Instituto Superior T\'ecnico, Universidade de Lisboa \newline \indent Av. Rovisco Pais, 1049-001 Lisboa, Portugal}}
\email{htavares@math.ist.utl.pt}

\author[Pedro J.\ Freitas]{Pedro J.\ Freitas}\thanks{}
\address{\footnotesize{Pedro J.\ Freitas \newline \indent CEAFEL and Department of Mathematics \newline \indent Faculdade de Ci\^encias, Universidade de Lisboa}}
\email{pjfreitas@fc.ul.pt}


\begin{abstract}
There are known constructions for some regular polygons, using only straightedge and compass, but not for all polygons---the Gauss-Wantzel Theorem states precisely which ones can be constructed.

The constructions differ greatly from one polygon to the other. There are, however, general processes for determining the side of the $n$-gon (approximately, but sometimes with great precision), which we describe in this paper. We present a joint mathematical analysis of the so-called Bion and Tempier approximation methods, comparing the errors and trying to explain why these constructions would work at all.
\end{abstract}
\date{\today}
\maketitle

\section{Geometric constructions for the regular polygons}

As many of us have learned in high school, there are geometric constructions for some regular polygons, using only straightedge and compass. This is a classical problem, which is interesting on its own right, even though the methods have become obsolete with geometry computer programs. 

The constructions for the triangle, the square and the hexagon are simple. For the pentagon, the construction is more delicate, but well known. For the octagon, one can simply bisect the angle of the square. The same goes for the decagon, one bisects the angle of the pentagon.

The reader may have noticed that we have skipped two regular polygons: the heptagon (7 sides) and the nonagon or enneagon (9 sides). There is a good reason for it: no matter how hard we try, we will never find a straightedge-and-compass construction for them. This is a consequence of the following classical result.\medskip

{\bf Gauss-Wantzel's Theorem}. {\it The division of the circle in $n$ equal parts with straightedge and compass is possible if and only if 
$$n=2^k p_1\ldots p_t$$
where $p_1,\ldots ,p_t$ are distinct Fermat primes.}\medskip

A {\it Fermat prime}\ is a prime of the form $2^{2^m}+1$. Presently, the only known Fermat primes are 3, 5, 17, 257 and 65537.\medskip

Gauss proved, in his early years, that the 17-gon is constructible. He went on to formulate the theorem, Wantzel concluded the proof in \cite{Wantzel}.\footnote{See \cite{hepta} for a video of David Eisenbud doing the construction of the 17-gon and \cite{hepta2} for a discussion of the mathematics involved.}  We refer to chapter 19 of the book by Stewart \cite{Stewart} for a complete proof of this theorem, where you can also find a reference for the construction of the 257-gon (!), as well as some funny anecdotes about the 65537-gon.\medskip

As we said, this result implies that the heptagon and the nonagon are not constructible with straightedge and compass:  7 is not a Fermat prime, and $9=3\times 3$, with the Fermat prime 3 appearing twice in the factorisation of 9. Nevertheless, it is possible to find good approximate constructions for both these polygons.\medskip

\begin{figure}[htb]
\centerline{\includegraphics[scale=0.70]{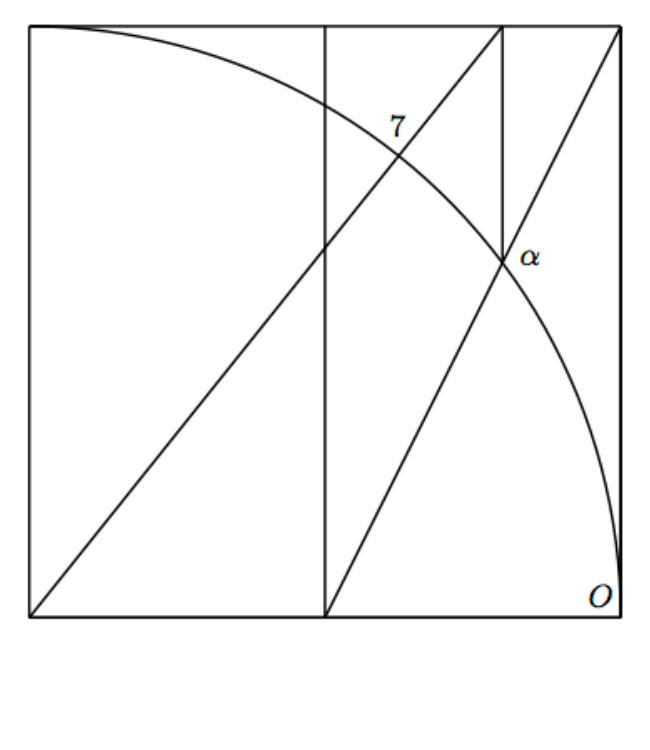} \qquad 
\includegraphics[scale=0.72]{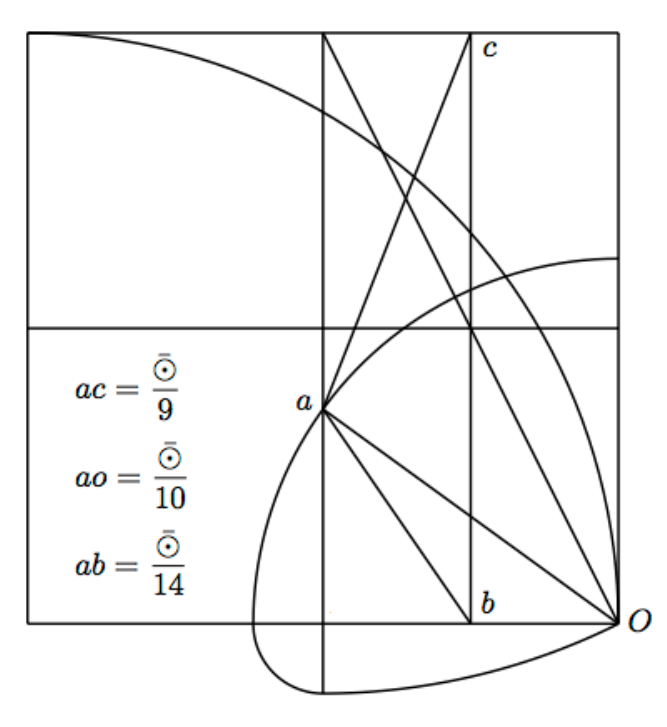}} 
\caption{Reproduction of drawings of Almada Negreiros}
\label{Almada}
\end{figure}

Several artists took an interest in finding ways to divide the circle in $n$ parts, and use this in their work. D\"urer was known for his taste for these constructions (see \cite{durer}). In the 20th century, there was a famous Portuguese modernist artist, Almada Negreiros, who produced drawings consisting solely of such constructions. We reproduce two of them here in Figure \ref{Almada} (the originals can be found in \cite{Almada1,Almada2}).

In the first one, the arc $\bigfrown{O7}$ is the 7th part of the circle, and in the second, the lines $ab$, $aO$ and $ac$ are the chords for the 14th, 10th and 9th parts of the circle, respectively. The errors are $1\%$ for the 14th part, $0.2\%$ for the 7th and an amazing $0.001\%$ for the 9th. The 10th part is exact. See \cite{Fr} and \cite{FrSi2} for a detailed analysis.

\section{A general construction}

\begin{figure}[htb]
\centerline{\includegraphics[scale=0.6]{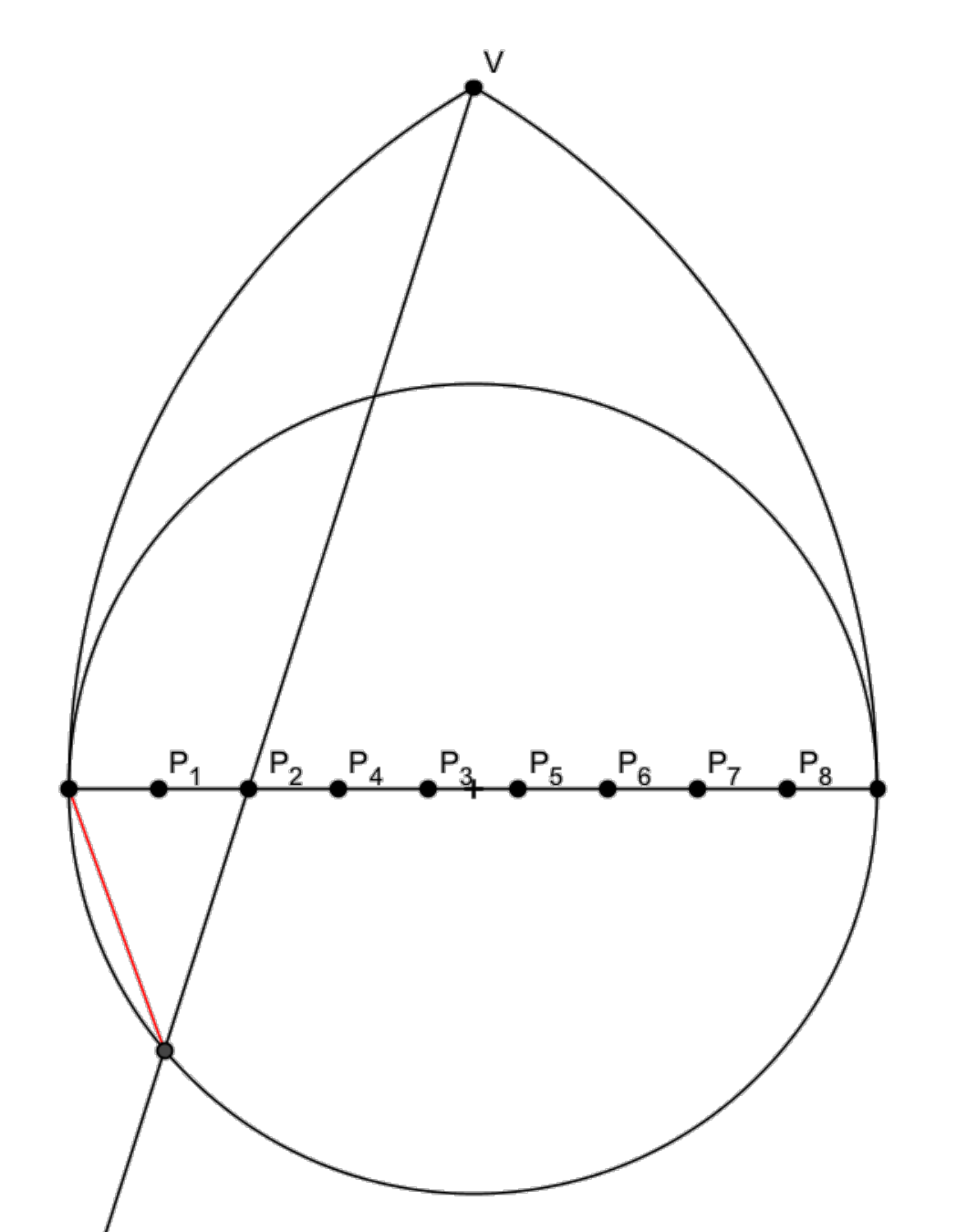}}
\caption{Bion method for the nonagon}
\label{bion_nonagon}
\end{figure}

Generally, the constructions become more and more intricate as we increase the number of sides of the polygon. Moreover, they are very distinct from one polygon to the other. While searching for good approximations, one must try to find a compromise between the complexity of the construction and its accuracy. As it happened, one of the authors once asked a friend who teaches the high school descriptive geometry class if there was a general approximate construction for the $n$-gon which such characteristics. And it turns out there was. She proceeded to draw the picture in Figure \ref{bion_nonagon} as a construction for the nonagon.\medskip

\begin{figure}[htb]
\centerline{\includegraphics[scale=0.6]{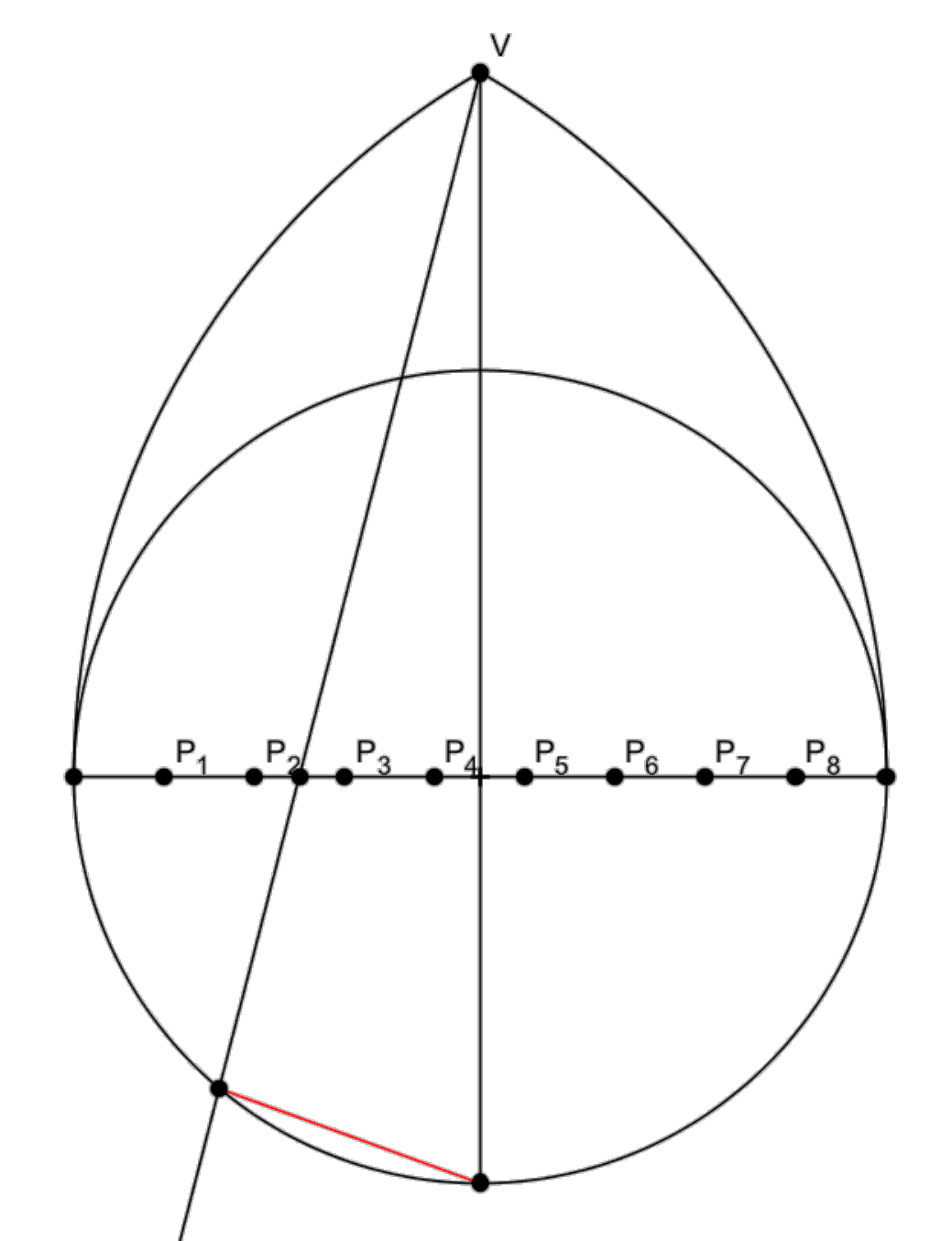}}
\caption{Tempier method for the nonagon}
\label{tempier_nonagon}
\end{figure}

After drawing a circle, one determines the intersection point $V$ of two arcs of circle centred at the endpoints of a diameter, having this diameter as radius.\footnote{This figure is usually called the {\em vesica  piscis}, the fish's bladder, and appears very frequently in the composition of Medieval art, for instance.} Then one divides the diameter in 9 parts and takes the {\em second point}\ Êfrom the left to define the ray. The process is completely general --- for the $n$-gon one simply has to divide the diameter in $n$ parts an always take the second point from the left.

This method can be found in some descriptive geometry books, such as \cite{brasileiro} and \cite{antigo}: it is called the {\em Bion method}, cf.\ the book \cite{Bi} and \cite{Ho}. A variation of this method was proposed by Tempier \cite{Te1,Te2}: one can also take a different ray, directed by another point, situated also on the diameter at a distance of two $n$-th parts from the centre of the circle. This is called the {\em Tempier method}, which we exemplify for the nonagon in Figure \ref{tempier_nonagon}.

The 19th century papers \cite{Ho,Te1, Te2} include mathematical descriptions of the methods and error tables, and notice that for the 17-gon, the construction is much simpler that Gauss's (even though it is approximate). The subject seems to be not so well known in the present mathematical community, and these constructions have received recent interest in the article \cite{Mi}, where the mathematics involved is again described.

For a visual approach to the methods, we refer the reader to the following online worksheets, that show the polygons determined by Bion's method, \cite{Bion_ggb}, and Tempier's method, \cite{Tempier_ggb}.\medskip

Still, we were puzzled about why would these methods afford an approximation of the $n$-th part of the circle at all? What was so special about the point $V$ and the length of 2 $n$-th parts of the diameter?

Therefore, our main aim, for the remainder of the paper, is to present a unified mathematical analysis of these methods (slightly different from the ones presented in previous papers) in order to provide the answers to these two questions.

\section{Measuring the error} 

From now on, we will consider the errors in the angles determined by the constructions. To know these errors, we must first find the exact values of the angles determined by each one. The analysis presented here is slightly different from the one in the original papers, and is thought so that it can be used for both methods.

\begin{figure}[htb]
\centerline{\includegraphics[scale=0.65]{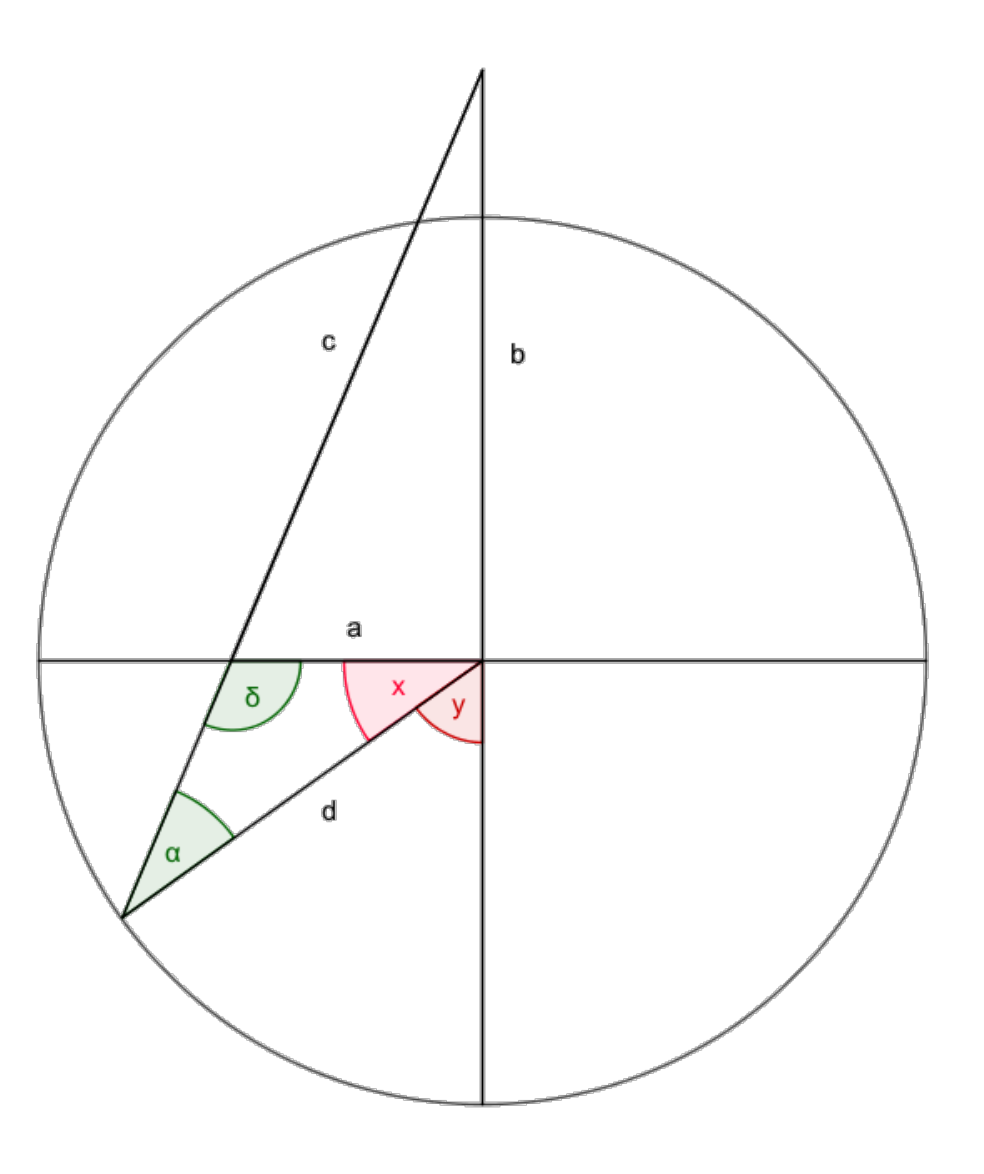}}
\caption{Measuring angles $x$ and $y$}
\label{generalfigure}
\end{figure}

 Consider Figure \ref{generalfigure}, where the vertical and horizontal lines meet at the centre.

We want the values of angles $x$ and $y$. We start studying the angle $\delta$:
\[
\sin \delta  =  \sin (\pi -\delta) = \frac bc, \qquad\cos \delta  =  -\cos (\pi -\delta) = -\frac ac.
\]
Since $\delta$ is an obtuse angle, given the main restrictions of sine, cosine and tangent we have
\begin{equation}
\label{delta}
\delta = \pi-\arcsin \frac bc = \arccos \left(-\frac ac\right) = \pi+\arctan \left(-\frac ba\right).
\end{equation}
This identifies the angle $\delta$ from the lengths $a$, $b$, and $c$. Now we establish the value of angle $\alpha$, using the law of sines.
$$\frac{\sin \alpha}a = \frac{\sin \delta} d \eq \sin \alpha = \frac{a\sin \delta} d \eq \sin \alpha = \frac{ab}{cd}.$$
So
\begin{equation}
\label{alpha}
\alpha =\arcsin \frac{ab}{cd}.
\end{equation}
By considering the figure, we see that 
$$x=\pi-\delta-\alpha\quad \text{and} \quad y=\frac \pi 2 -x = \delta + \alpha -\frac\pi 2.$$
Now we can calculate the value of $x$: recalling the expressions \eqref{delta} and \eqref{alpha},  
\begin{equation}
\label{angulo_x}
x = \pi - \left(\pi - \arcsin \frac bc\right) - \arcsin \frac{ab}{cd} = \arcsin \frac bc -\arcsin \frac{ab}{cd}.
\end{equation}
For the value of $y$, notice that
\[
\cos\left(\frac{\pi}{2}-\alpha\right)=\sin\alpha=\frac{ab}{cd},
\]
so, since $\frac \pi2-\alpha$ is an angle in $[0,\pi]$, 
\begin{equation}
\label{angulo_y}
y = \delta + \alpha -\frac\pi 2 = \arccos \left(-\frac ac\right) - \arccos \frac {ab}{cd}.
\end{equation}

We can now express the values of the angles that appear in the Bion and Tempier methods: it suffices to calculate the values of $a$, $b$, $c$ and $d$. 

Let's take as unit the radius of the circle, i.e., $d=1$, since it will be clear from the computations that the angles do not depend on this value. 

\subsection*{Bion method.} In this case, we are interested in the value of $x$. We have $a=1-4/n=(n-4)/n$ and $b=\sqrt{2^2-1^2}=\sqrt{3}$, so that
\[
c=\sqrt{(\sqrt{3})^2+\left(\frac{n-4}{n}\right)^2}=\frac{2\sqrt{n^2-2n+4}}{n}.
\]
Therefore,
\[
x=\arcsin\frac{\sqrt{3}n}{2\sqrt{n^2-2n+4}}-\arcsin \frac{\sqrt{3}(n-4)}{2\sqrt{n^2-2n+4}}
\]
We present in the table below the values, for $n=4,\ldots, 20$, of $2\pi/n$, the approximation $x$, as well as the error and the relative error. All the quantities were rounded to four decimal digits.

{\small \begin{center}
\begin{tabular}{ c | c | c | c | c }
  & Angle & Approximation & Error & Relative Error \\
   $n$    &	   $\frac{2\pi}{n}$ &		$x$          & $\frac{2\pi}{n}-x$ & $\frac{|2\pi/n-x|}{2\pi/n}$\\     \hline \hline
 4    & $1.571$	  & 1.571		&  0				&  0 	\\ \hline
5    &	 $1.257$          & 1.256	& 0.0008	&     0.0006\\\hline
6 &	$1.047$ & 1.047   &  0  			&	0\\\hline
7 &	$0.8976$ & 0.8992    &-0.0016		& 0.0017\\\hline
8 & $0.7854$	 & 0.7887	& -0.0033		& 0.0042\\\hline
9 & $0.6981$	 & 0.7030	&  -0.0048	         &0.0069\\\hline
10 & $0.6283$	& 0.6345    &	   -0.0062	&0.0099\\\hline
11 & $0.5712$ 	& 0.5785	&  -0.0073		&0.0129\\\hline
12 & $0.5236$	& 0.5319     &-0.0083		&0.0158\\\hline
13 & $0.4833$	& 0.4923   &  -0.009		&0.0186\\\hline
14 & $0.4488$	& 0.4584  &	-0.0096	&0.0214\\\hline
15 & $0.4189$	& 0.4289	& -0.01		&0.024\\\hline
16 & $0.3927$	& 0.4031	& -0.0104		&0.0265\\\hline
17 & $0.3696$	& 0.3803	& -0.0107		&0.0288\\\hline
18 & $0.3491$	& 0.3599	& -0.0108		&0.0311\\\hline
19 & $0.3307$	& 0.3417	& -0.011		&0.0332\\\hline
20 & $0.3142$	& 0.3252	&-0.0111		&0.0352
\end{tabular}
\end{center}}

Elementary --- but rather tedious --- computations show that the function ``relative error'', defined on the integers greater than or equal to 4, given by
\[
n\mapsto \frac{\frac{2\pi}{n}-x(n)}{\frac{2\pi}{n}}=1-\frac{nx(n)}{2\pi}
\]
 is only exact for $n=4,6$, has a maximum for $n=5$, and it is strictly decreasing for $n\geq 6$, converging as $n\to \infty$ to $1-2\sqrt 3/\pi\simeq -0.1026$. Thus Bion's formula is only exact for $\pi/2$ and $\pi/3$, and in the other terms the error can go from $0.1\%$ to $10.3\%$ (although, as we can see from the table above, it is still reasonable for $n=20$: approximately $3.5\%$). 

\subsection*{Tempier method.} In this method, $a=4/n$, $b=\sqrt{3}$ and 
$c=\sqrt{3n^2+16}/n$. In this case the approximation for $2\pi/n$ is given by
\[
y=\arccos \left(-\frac{4}{\sqrt{3n^2+16}}\right)-\arccos\frac{4\sqrt{3}}{\sqrt{3n^2+16}}.
\]
{\small
\begin{center}
\begin{tabular}{ c | c | c | c | c }
   & Angle & Approximation & Error & Relative Error \\
    $n$  	&$\frac{2\pi}{n}$ &		$y$          & $\frac{2\pi}{n}-y$ & $\frac{|2\pi/n-y|}{2\pi/n}$\\     \hline \hline
4 & $1.571$	&1.571	&0  &	0 	\\ \hline
5 &	$1.257$	&1.246	&0.0111	&0.0088        \\\hline
6 &	$1.047$	& 1.039	&0.0083	&0.0079\\\hline
7 &	$0.8976$ & 0.8923	&0.0053	&0.0059\\\hline
8 & $0.7854$ & 0.7821	&0.0033	&0.0042\\\hline
9 & $0.6981$	 & 0.6962	&0.0019	&0.0027\\\hline
10 & $0.6283$	& 0.6273	&0.001	&0.0016\\\hline
11 & $0.5712$	& 0.5708	&0.0004	&0.0007\\\hline
12 & $0.5236$	& 0.5236	&0	&0\\\hline
13 & $0.4833$	& 0.4836	&-0.0003	&0.0006\\\hline
14 & $0.4488$	& 0.4493	&-0.0005	&0.001\\\hline
15 & $0.4189$	& 0.4195	&-0.0006	&0.0014\\\hline
16 & $0.3927$	& 0.3934	&-0.0007	&0.0017\\\hline
17 & $0.3696$	& 0.3703	&-0.0007	&0.002\\\hline
18 & $0.3491$	& 0.3498	&-0.0008	&0.0022\\\hline
19 & $0.3307$	& 0.3315	&-0.0008	& 0.0024\\\hline
20 & $0.3142$	& 0.315	&-0.0008 &	0.0026
\end{tabular}
\end{center}}

By studying now the function 
\[
n\mapsto \frac{\frac{2\pi}{n}-y(n)}{\frac{2\pi}{n}}=1-\frac{ny(n)}{2\pi},
\]
we see that it reaches a maximum at $n=5$, and then it is strictly decreasing for $n\geq 5$, converging to $-(6+2\sqrt{3}-3\pi)/(3\pi)\simeq -0.00417$. Then the formula is only exact for $\pi/4$ and $\pi/6$, while in the other terms the error is never worse that $0.9\%$! \medskip

Comparing the two methods, one concludes that Bion's is better for $n=5,6,7$, they are more or less equivalent for $n=8$, and then for $n\geq 9$ Tempier's is much more efficient.
\medbreak

\begin{figure}[htb]
\centerline{\includegraphics[scale=0.9]{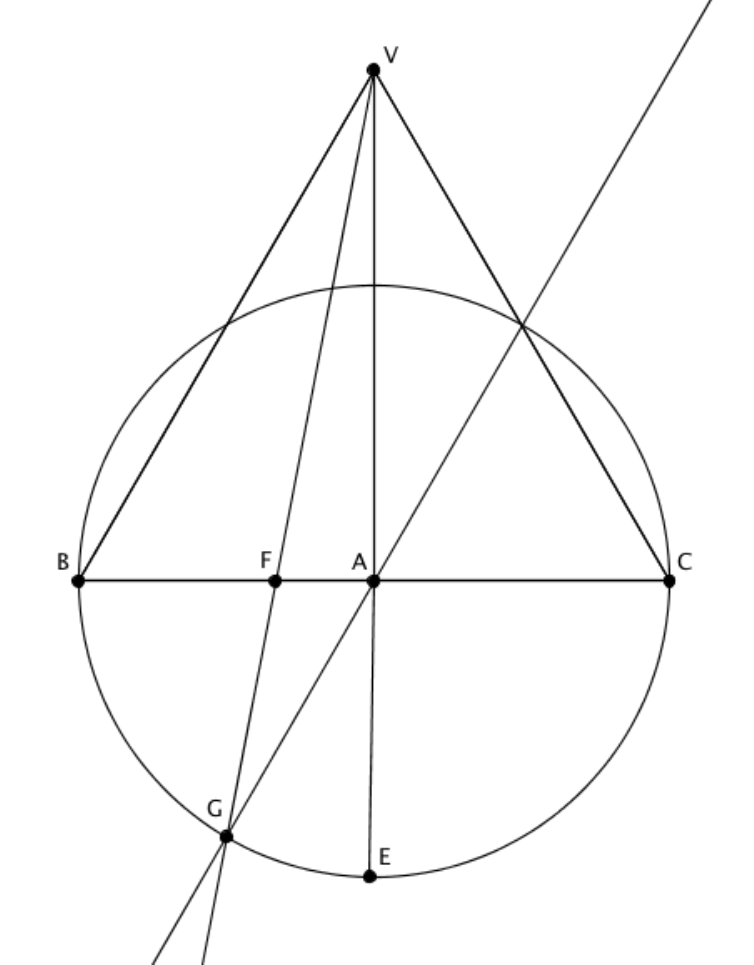}}
\caption{Exact cases}
\label{exactcases}
\end{figure}

We noticed that for the Bion hexagon and the Tempier dodecagon the construction was exact. Both polygons give origin to Figure \ref{exactcases}, with 
$$BF=\frac 13 BC,\quad FA=\frac 16 BC.$$
This also means that $BF=2FA$.

It is straightforward to prove this exactness using more elementary methods. Take a parallel to $BV$ through point $A$, and let $G$ be its intersection with line $VF$. Using the similarity of triangles $BFV$ and $GFA$, one can show that $\measuredangle BAG=\pi/3$ (which implies that $\measuredangle GAF=\pi/6$), and that point $G$ is on the circle.

This proves what we want, since $\pi/3$ and $\pi/6$ are the centre angles for the hexagon and dodecagon, respectively.

\section{The rectification of the quadrant}

We now turn to the question regarding point $V$. Why was this point chosen for these constructions?\medskip

\begin{figure}[htb]
\centerline{\includegraphics[scale=0.8]{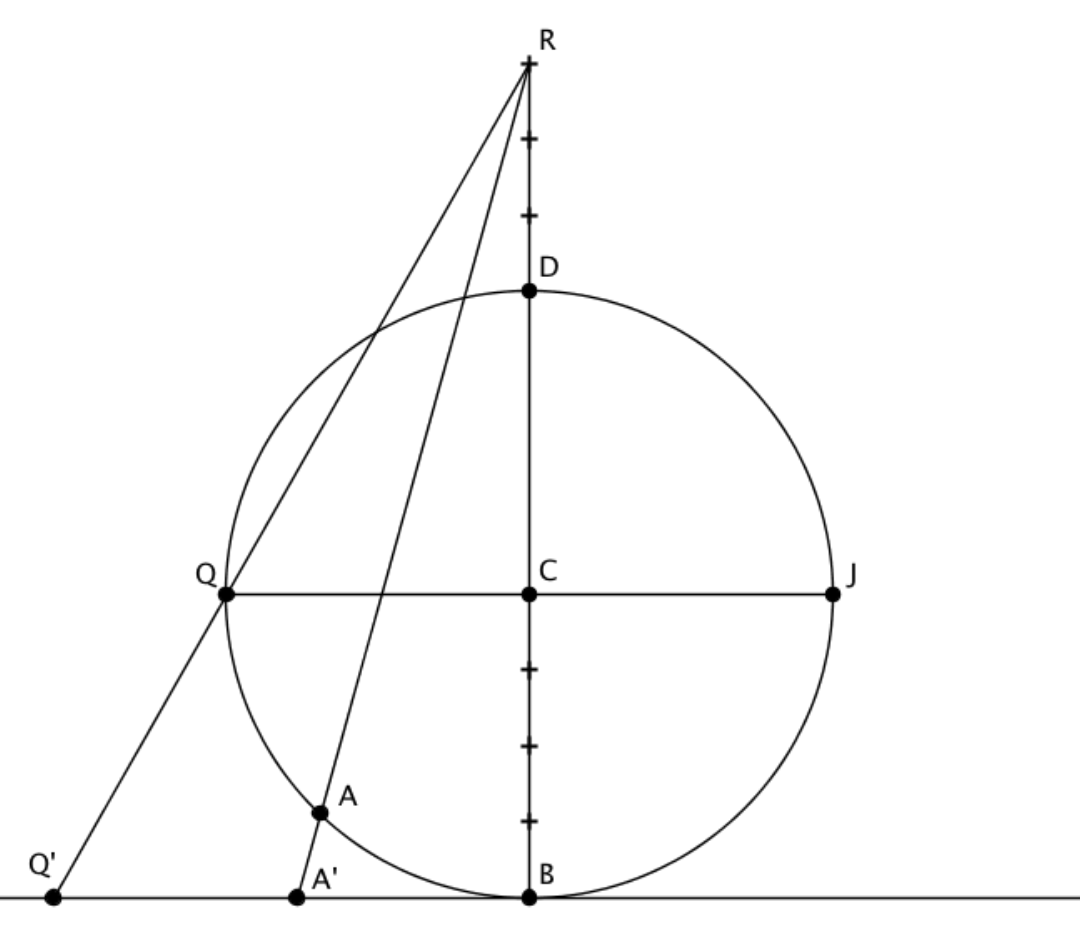}}
\caption{Rectification of arcs}
\label{quadrant}
\end{figure}

We noticed that, in \cite{brasileiro}, there is a construction for rectification of arcs smaller that $\pi/2$ that uses a very close point, see Figure \ref{quadrant}. The point $R$ is marked on the vertical line at a distance from $D$ equal to 3/4 of the radius, as the figure shows (we used $R$ for rational).\medskip

This construction is used to produce a segment $A'B$ with length approximately equal to that of the arc $\bigfrown{AB}$. Without delving too much on the accuracy of this construction, we just calculate what is the length obtained when the arc is a full quadrant, $\bigfrown{QB}.$ 

By simple proportions, and still considering the circle to have unit radius, we have
$$\frac{Q'B}{RB} = \frac{QC}{RC} \eq Q' B = \frac{11/4}{7/4}Ê\eq Q'B = \frac{22/7} 2.$$ 
So, the position of point $R$ is such that the approximation of $\pi$ for the rectification of the whole quadrant is the famous $22/7\simeq3,14286$. \medskip

If we use the point $V$ from the previous constructions to rectify the quadrant, and do the corresponding proportion, we find that the implicit approximation of $\pi$ is $(2+2\sqrt 3)/3 \simeq 3,15470$.\medskip

One can also work the other way around: what would be the position of a point $P$, on the line $BD$, such that the rectification of the full quadrant would be exact? If we work out the proportion once more, we find that the distance from $P$ to $C$ would have to be $2/(\pi-2) \simeq 1,75194$. This is a non-constructible length, since it is transcendental, and point $R$ is a very good approximation. \medskip

We summarise our computations in the following table.

{\renewcommand{\arraystretch}{1.3}
$$\begin{array}{c|c|c}
\text{Point} & \text{Distance to } C & 2 \times \text{(rectified quadrant)}\\ \hline
V & \sqrt 3 \simeq 1,73205 & (2+2\sqrt 3)/3\simeq 3,15470\\ \hline
R & 7/4 = 1,75 & 22/7 \simeq 3,14286\\ \hline
P & 2/(\pi-2) \simeq 1,75194 & \pi \simeq 3,14159
\end{array}$$
}

If we use this point $P$ for the Tempier process, and calculate the corresponding relative error, we see that the limit is 0 as $n\to \infty$.
\medskip

We conclude that all points $V$, $R$ and $P$ are reasonable base points for rectifying arcs.  In particular, it probably explains why $V$ was chosen: of the constructible lenghts, it is probably the easiest to mark. 

So, when we look for the value of the angle in the Tempier  or Bion methods, we can look at the length of the corresponding (approximately) rectified segment on the tangent line. 

Now, if you divide the diameter into $n$ equal parts, this corresponds to a division of the full rectified half-circle into $n$ equal parts as well. Each part will then have length $\pi/n$, so one has to take two of these in order to get a length of $2\pi/n$, as we wish. And this explains why we must take 2 segments on the diameter in order to define the angle $2\pi/n$.

\section{Conclusion}

In this paper we looked at geometric constructions that are familiar to artists, architects and people interested in descriptive geometry. These constructions don't seem to be so well known in the mathematical community. We studied the mathematics involved in these constructions, reformulating the mathematics that appears in the original papers, in order to understand how they work.

In practice, approximations are acceptable in certain circumstances, and some errors are so small that the very thickness of the trace renders them negligible. Our computations confirm that the Bion and Tempier methods are good options for the $n$-gon, being acceptable even in the cases when an exact construction is known.

\end{document}